\documentclass[12 pt]{article}
\usepackage{graphicx} 
\usepackage{amssymb,amsthm,amsmath}
\usepackage{mathrsfs}
\usepackage{mathtools}
\usepackage{fullpage}
\usepackage{microtype}
\usepackage{indentfirst}  
\usepackage{titlesec}
\usepackage{etoolbox}
\usepackage{hyperref}
\usepackage{fancyhdr} 

\newtheorem{theorem}{Theorem}[]
\numberwithin{theorem}{section}

\newtheorem{lemma}[theorem]{Lemma}

\newtheorem{corollary}[theorem]{Corollary}

\newtheorem*{theorem*}{Theorem}
\newtheorem{remark}[theorem]{Remark}

\newcommand{\Fq}{\mathbb{F}_q}
\newcommand{\R}{\mathbb{R}}
\newcommand{\Q}{\mathbb{Q}}
\newcommand{\Z}{\mathbb{Z}}

\newcommand{\Orb}{\mathcal{O}}
\newcommand{\tf}{transverse-free}

\newcommand{\AG}{\operatorname{AG}}
\newcommand{\PG}{\operatorname{PG}}
\newcommand{\Proj}{\operatorname{Proj}}
\newcommand{\Pp}{\mathbb{P}}
\newcommand{\A}{\mathbb{A}}

\newcommand{\homog}{\mathrm{hom}}

\newcommand{\Sing}{\mathrm{Sing}}

\titleformat{\section}
  {\normalfont\large\centering} 
  {\normalfont\thesection} 
  {1em} 
  {}
\titleformat{\subsection}
  {\normalfont\large} 
  {\normalfont\thesubsection} 
  {1em} 
  {}
\providecommand{\keywords}[1]
{
  \small	
  \textbf{\textit{Keywords---}} #1
}

\title{
\textsc{Almost All Transverse-Free Plane Curves Are Trivially Transverse-Free}
}
\author{
  Alejandro Lopez\thanks{Rice University, \href{mailto:atl5@rice.edu}{atl5@rice.edu}}
  \and Bella Villarreal\thanks{Grinnell College, \href{mailto:villarre2@grinnell.edu}{villarre2@grinnell.edu}}
  \and Ren Watson\thanks{University of Texas at Austin, \href{mailto:renwatson@utexas.edu}{renwatson@utexas.edu}}
  \and Jaedon Whyte\thanks{Massachusetts Institute of Technology, \href{mailto:jnw@alum.mit.edu}{jnw@alum.mit.edu}}
}
\begin{document}
\maketitle
\linespread{1.03}
\begin{abstract}

   Call a curve $C \subset \mathbb{P}^2$ defined over $\mathbb{F}_q$ \textit{transverse-free} if every line over $\mathbb{F}_q$ intersects $C$ at some closed point with multiplicity at least 2. In 2004, Poonen used a notion of density to treat Bertini Theorems over finite fields. In this paper we develop methods for density computation and apply them to estimate the density of the set of polynomials defining \textit{transverse-free} curves. In order to do so, we use a combinatorial approach based on blocking sets of $\PG(2, q)$ and prove an upper bound on the number of such sets of fixed size $< 2q$. We thus obtain that nearly all transverse-free curves contain singularities at every $\Fq$-point of some line.
   
\end{abstract}
\keywords{transverse-free, blocking sets, finite fields, projective plane} 

MSC[2020] 11G20; 51E20; 51E21

\section{Introduction}

For a quasi-projective subscheme $X \subseteq \Pp^n$ over a field $k$, Bertini's Theorem states that if $X$ is smooth, then for a generic hyperplane $H \subset \Pp^n$ over $\bar{k}$, the intersection $H\cap X$ is also smooth. Over infinite fields $k$, we can find many such hyperplanes that are defined over $k$, but cannot conclude on their existence over finite fields. Poonen's work, \cite{poonen}, remedies this situation by showing that instead of considering generic properties of arbitrary hyperplanes, one can instead consider properties of hypersurfaces, $H_f = \{f = 0\}$, over $k=\Fq$. In this manner, Poonen proves a Bertini-type theorem for finite fields, which explicitly computes the fraction of homogeneous polynomials $f$ over $\Fq$ such that $H_f \cap X$ is smooth of dimension $\dim{X} - 1$. In the opposite direction, an Anti-Bertini theorem is also proved stating that for any finite collection of hypersurfaces $\{H_{f_i}\}$ over $\Fq$, there exist smooth hypersurfaces $X$ over $\Fq$ that intersect none of them smoothly with dimension $n-2$. The proof again argues that some positive fraction of defining polynomials are satisfactory. Thus the problem of computing this associated \textit{density} for some special classes of hypersurfaces over $\Fq$ becomes of interest. 

In this paper we treat the problem for \textbf{transverse-free} curves over $\Fq$ in $\Pp^2$, that is, one-dimensional hypersurfaces for which every line over $\Fq$ fails to have a smooth dimension $0$ intersection. This leads to the following question: what is the density of the set of polynomials defining transverse-free curves? This would extend the work of \cite{freidin-asgarli}, wherein the authors introduced and considered the case of transverse-free curves that are additionally smooth. From our initial perspective, the question can also be viewed as asking for a more precise analog of Poonen's Anti-Bertini Theorem, details of which may be found in the subsequent remark.

\begin{remark}
    In fact, Bertini's Theorem may be strengthened to show that geometrically integral transverse-free curves do not exist over infinite fields as the number of singularities must be finite. Using Proposition $2.7$ of \cite{poonen} stating that the density of geometrically integral curves is $1$, we therefore see that this is equivalent to asking for the density of counterexamples to this strengthened version of Bertini's Theorem for projective plane curves over finite fields. 
\end{remark} 

Henceforth we shall assume all subschemes are over our base field $\Fq$ and refer to $\Fq$-rational points simply as rational points. In order to state the question precisely, set $R = \Fq[x_0, x_1, x_2]$ and let $R_d \subseteq R$ be the linear subspace of homogeneous degree $d$ polynomials, and $R_{\homog} \coloneqq \bigcup_{d=1}^{\infty} R_d$. Then the \textbf{density} of a subset $\mathcal{P}\subseteq R_{\homog}$ is given by: \[\mu(\mathcal{P}) \coloneqq \lim_{d\to \infty} \frac{|\mathcal{P}\cap R_d|}{|R_d|}.\] We denote the set in question by
\[\Omega \coloneqq \{f\in R_{\homog} | C_f = \Proj(R/(f)) \subseteq \Pp^2 \text{ is transverse-free}\},\] so that we ask for the value of $\mu(\Omega)$.

Here and throughout, a \textbf{curve} means a closed subscheme of $\Pp^2$ of the form $C_f = \Proj(R/(f))$ for a nonzero $f\in R_{\homog}$. In Section \ref{sec:density} we will consider more generally \textbf{hypersurfaces} $H_f\subseteq \Pp^n$, and everything may be analogously defined so long as we instead consider the coordinate ring $R =\Fq[x_0, x_1, \dots, x_n]$ depending on the dimension $n$.

Notice that $H_{f_1} = H_{f_2}$ only when there exists $u \in \Fq^\times$ so that $f_1 = uf_2$ and therefore the maps $f \in R_d \setminus \{0\} \mapsto H_f$ are $(q-1)$-coverings. Considering that we have defined all hypersurfaces as arising from such a construction, it is in light of this fact that we are justified in speaking of the densities of various classes of hypersurfaces as the densities of the subsets of homogeneous polynomials that define them. Our main result is the following asymptotic formula for $\mu(\Omega)$ as the prime power $q$ grows large.

\begin{theorem}\label{theorem:TBS1}
    \[\mu(\Omega) = (q^2 + q+1)q^{-3(q+1)}\Bigl(1+O\bigl(q^{-3\sqrt{q} + 2}\bigr)\Bigr).\]
 In the simplest cases of $q = p$ and $q = p^2$ for a prime $p$, we find more precise estimates. In particular, \begin{align*}\mu(\Omega) &= (q^2 + q+1)q^{-3(q+1)}\left(1+O\left(\left(\frac{q}{2}\right)^{-\frac{1}{2}(q+1)}\right)\right), \allowdisplaybreaks \text{ and} \\ \mu(\Omega) &= \left[(q^2 + q+1)q^{-3(q+1)}+(q^3 + q\sqrt{q})(q+1)q^{-3(q+\sqrt{q}+1)}\right]\allowdisplaybreaks\left(1+O\left(\left(\frac{q}{2}\right)^{-\frac{1}{2}(q+1)}\right)\right),\end{align*} respectively.
\end{theorem}

Beyond the numerical statement of Theorem \ref{theorem:TBS1}, the proof also has interesting qualitative ramifications for transverse-free curves. Using density as a measure of frequency, we show that transverse-free curves are exceedingly rare, with their existence typically linked to combinatorial structures in the projective plane called \textbf{blocking sets}, that is, sets of rational points that intersect every line. This connection allows us to conclude that almost all transverse-free curves correspond to small blocking sets. In particular, the simplest blocking sets are referred to as \textbf{trivial blocking sets} and contain all of the rational points of some line. We call a curve \textbf{trivially transverse-free} if its set of  singular rational points is a trivial blocking set. The main estimate shows that for $q$ sufficiently large, nearly all transverse-free curves are trivially transverse-free as suggested by the paper title.

A key ingredient in the proof of Theorem \ref{theorem:TBS1} is a new result on the number of blocking sets of projective planes of certain sizes. We call a blocking set minimal if no proper subset of it is a blocking set. Then we have the following result.
\begin{theorem}\label{theorem:min-blset-bound} Let $B_k$ denote the number of \textbf{minimal blocking sets} of size $k$. For $q \leq k \leq 2q$, 
\begin{equation*}
     B_k \leq \frac{\binom{q^2 + q + 1}{2k - 2q}}{\binom{k}{2k-2q}}.
\end{equation*}
\end{theorem}
In the $q=p$ and $q=p^2$ cases, we have stricter bounds from the literature on $B_k$ and correspondingly more precise estimates in Theorem \ref{theorem:TBS1}.
 
\subsection{Paper Outline}

In Section \ref{sec:density} we address the question of computing $\mu$ for classes of hypersurfaces determined by first-order local constraints in general. Starting from the basic property of transversality to a smooth subscheme, we define a family of subsets of $R_\homog$ for which $\mu$ is explicitly computable. We will call these the simple sets. We also introduce a notion of conditional density which will provide a convenient formulation for the material of Section \ref{sec:Sam-Sung’s-Note}. In Subsection \ref{subsec: rationality}, we develop more sophisticated notions concerning these subsets and use them to prove Theorem \ref{thm:inQ} showing that in the dimension $2$ case, the density of all subsets in the Boolean algebra generated by simple sets can be computed using the introduced material. 

In Section \ref{sec:TTFC} we exhibit foundational applications of the methods in the previous section in the investigation of trivially transverse-free curves. The ensuing heuristics will guide the proof of the main theorem. In Section \ref{sec:Sam-Sung’s-Note} we provide our major tool for estimating density, Theorem \ref{thm:Sam-Sung’s-Note}. In particular, we prove a sub-independence result for transverse conditions under some smoothness assumption and give some applications. These include Theorem \ref{thm:smBounds}, which improves on an upper bound for the density of smooth transverse-free curves first exhibited in \cite{freidin-asgarli}. 

In Section \ref{sec:block-time}, we collect results from the combinatorial literature on blocking sets. In particular, we look for enumerative results that lend themselves to computing the conditional density of $\Omega$ with respect to those curves where the set of rational singular points is prescribed. As Section \ref{sec:TTFC} suggests, this leads naturally to the question of how many lines avoid all of these singularities. As already referenced, the new work in this section consists of the proof of Theorem \ref{theorem:min-blset-bound}, which effectively limits the contribution of most ways of arranging the rational singular points to $\mu(\Omega)$.

Finally, Section \ref{sec:main-proof} combines the results of the previous two sections to prove Theorem \ref{theorem:TBS1} by direct computation.

\section*{Acknowledgements}
This work is supported by the National Science Foundation under Grant Number DMS-1851842, and conducted at the MathILy-EST 2023 REU under the supervision of Dr. Brian Freidin. Any opinions, findings, and conclusions or recommendations expressed in this paper are those of the authors and do not necessarily reflect the views of the National Science Foundation.
The authors declare no conflicts of interest. 

\section{ 
General Methods for Density Computation}\label{sec:density}

As the name transverse-free suggests, we understand the subscheme properties we are interested in as stemming from the notion of a transverse intersection of subschemes. There does not appear to be a standardized definition for this, so we shall say that a subscheme, $H \subseteq \Pp^n$, is transverse to a smooth subscheme, $X \subseteq \Pp^n$, when the scheme-theoretic intersection, $H \cap X$, is a smooth subscheme of $X$ with codimension equal to the codimension of $H$. It follows that a curve is transverse-free exactly when there are no transverse lines. 

We will only investigate the case where $H$ has codimension $1$, as this is the setting in which we have a simple notion of density which is used in \cite{poonen}. We will employ the conventions from the same paper that the empty scheme, $\emptyset$, is smooth of every dimension and a subscheme $X$ is smooth of every dimension at any closed point $P \notin X$. Smoothness and dimension are local properties, so we can write the set of $f$ such that $H_f$ is transverse to $X$ as \[\{f \in R_\homog \colon H_f \, \cap \, X \text{ is smooth of dimension }\dim X - 1 \text{ at all } P \in X\}.\] We shall think of this as imposing local transverse intersection constraints on $H_f$ for each $P \in X$. 

This formulation inspires the possibility of only imposing a subset of these local constraints by considering a subscheme $U \subseteq X$. Notice additionally that when $X \subseteq X'$, a transverse intersection with $X$ at $P$ implies a transverse intersection with $X'$ at $P$. We therefore find it useful for technical reasons to consider a ``maximal subscheme" that intersects every subscheme transversely. To that end we will make special allowance for the case where $X = \emptyset$ despite the fact that for subschemes $X$ we always have $\emptyset \subseteq X$ but usually not $X \subseteq \emptyset$. 

So for pairs of subschemes $(U, X)$ of $\Pp^n$ satisfying either $U \subseteq X$ with $X$ smooth or $X = \emptyset$, we will define the following subset of $R_\homog$: \[S_U^X \coloneqq \{f \in R_\homog \colon H_f \cap X \text{ is smooth of dimension } \dim X - 1 \text{ at all } P \in U \}.\] We will call any subset of $R_\homog$ that can be constructed in this manner an \textbf{elementary set}. As suggested, we have the following basic properties regarding the behaviors exhibited in the upper and lower entries of the notation:
\begin{enumerate}
    \item $S_U^{\emptyset} = R_\homog$,
    \item if $\emptyset \neq X \subseteq Y$, then $S_U^X \subseteq S_U^Y$.
    \item if in (2.), $X \subseteq Y$ open, then $S_U^X = S_U^Y$, and
    \item if $V \subseteq U$ closed, then $S_U^X = S_V^X \cap S_{U - V}^X$
\end{enumerate}

It is clear by the preceding elaboration that elementary sets are examples of sets whose density was studied by Poonen in \cite{poonen}. We therefore have means of computation from the same paper. The basic idea is that local conditions specified at different points should act independently of each other. More precisely, we say that a collection of subsets of $R_\homog$, $\{S_i\}_{i \in I}$, are \textbf{independent} if \[\mu\left(\bigcap_{i \in I} S_i\right) = \prod_{i \in I}\mu(S_i).\] Theorem 1.3 of \cite{poonen} allows us to conclude the following with little trouble. 

\begin{corollary}[\cite{poonen}]\label{cor:elem-indep}
    Suppose $\{U_i\}_{i \in I}$ are mutually disjoint. Then $\{S_{U_i}^{X_i}\}_{i \in I}$ are independent.
\end{corollary}

For concrete density values, we can apply Theorem 1.2 of \cite{poonen} to obtain the density in terms of the Weil zeta function $\zeta_U(s)$. We will make use of the case where $U$ is a linear projective subspace of $\Pp^n$.

\begin{corollary}\label{cor:zeta}
    $\mu(S_U^X) \in \Q$ and in particular it can be expressed in terms of the Weil zeta function of $U$ as \[\mu(S_U^X) = \zeta_{U}(\dim{X} + 1)^{-1}.\] If $U \simeq \Pp^m$ this obtains the explicit formula \[\mu(S_U^X) =  \prod_{k=0}^m (1 - q^{k - \dim{X} - 1}).\]

    \begin{proof}
        Though not a claim made explicitly in \cite{poonen}, small modifications can be made to the proof of Theorem 1.1 from that paper so that it computes $S_U^X$ generally instead of merely $S_X^X$. We however give a self-contained account. Since $U$ is a subscheme of $X$, it is a closed subscheme of an open subscheme $Y \subseteq X$, and it therefore suffices to consider the case where $U \subseteq X$ is closed. Then applying Corollary \ref{cor:elem-indep}, Theorem 1.1 of \cite{poonen}, and the fact that the Weil Zeta function is multiplicative gives \[\mu(S_U^X) = \frac{\mu(S_X^X)}{\mu(S_{X - U}^X)} = \frac{\mu(S_X^X)}{\mu(S_{X - U}^{X - U})} = \frac{\zeta_X(\dim X + 1)^{-1}}{\zeta_{X - U}(\dim X + 1)^{-1}} = \zeta_U(\dim X + 1)^{-1}.\] as desired. Rationality is a result of \cite{dworkin} and the other properties of the Weil zeta function may be referred to from \cite{zeta}.
    \end{proof}
\end{corollary} 

We will denote subsets of $R_\homog$ that do \textit{not} give rise to hypersurfaces transverse to a given smooth subscheme $X$. In fact, we will provide a more general notation. A \textbf{complementary set} is any subset of the form \[S_U^{X, Y} = S_U^X \setminus S_U^Y \text{ for } X \supseteq Y \text{ or } X = \emptyset,\] and in particular $S_U^{\emptyset, X}$ is the complement of $S_U^X$ in $R_\homog$. Via the definition, it can be checked that $\mu$ satisfies additivity over finite collections of disjoint subsets of $R_\homog$. As such, our definition deliberately requires that $S_U^X \supseteq S_U^Y$ so that the density is given easily as $\mu(S_U^{X, Y}) = \mu(S_U^X) - \mu(S_U^Y)$. Collectively we will refer to both elementary and complementary sets as \textbf{simple sets}, and it will be convenient to refer to a simple set that can be written in the form $S_U^X$ or $S_U^{X, Y}$ as being \textbf{supported} at $U$. 

While the definition of the collection of simple sets may appear contrived, they are versatile enough to encode many properties of hypersurfaces that depend on first-order local behavior. So for example, in the $n = 2$ case, the possibilities for the first-order local behavior of the curve $C_f$ at a rational point, $P$, are all captured by simple sets, namely, either $P \notin C_f$, $P \in C_f$ is a singular point, or $P \in C_f$ and there is a unique ``tangent" line $\ell$ that will have a multiple intersection with $C_f$ at $P$. This is the content of the following lemma which is useful for Section \ref{sec:TTFC}.

\begin{lemma}\label{lem:numbers}
    Fix a rational point $P$ and a line $\ell$ through it. Then the following sets are simple and have densities computed by Corollary \ref{cor:zeta} and the subtraction rule as follows:
    \begin{enumerate}
        \item $\mu(\{f \in R_\homog : P \notin C_f \}) = \mu(S_P^P) = 1 - q^{-1}$,
        \item $\mu(\{f \in R_\homog : \ell \text{ is the unique line with a multiple intersection with } C_f \text{ at } P\})$ $= \mu(S_P^{\Pp^2, \ell}) = q^{-2} - q^{-3}$,
        \item $\mu(\{f \in R_\homog : C_f \text{ has a singularity at } P \text{ or } C_f = \Pp^2\}) = \mu(S_P^{\emptyset, \Pp^2}) = q^{-3}$.
    \end{enumerate}
\end{lemma}
Notice that we do in fact have a partition $R_\homog = S_P^P \sqcup \left(\bigsqcup_{\ell \ni P} S_P^{\Pp^2, \ell}\right) \sqcup S_P^{\emptyset, \Pp^2}$ and correspondingly $1 = (1 - q^{-1}) + (q + 1)(q^{-2} - q^{-3}) + q^{-3}$. That we have an analogous situation concerning the first-order local behavior of $H_f$ at any closed point $P$ is the main content of the proof of Lemma \ref{lem:splitting}.

We may also extend the independence result elementary sets to simple sets, though it is necessary to restrict the number of complementary sets allowed to be finite.

\begin{corollary}\label{cor:simple-indep}
    Suppose $\{U_i\}_{i \in I} \cup \{V_j\}_{j = 1}^m$ are mutually disjoint. Then $\{S_{U_i}^{X_i}\}_{i \in I} \cup \{S_{V_j}^{X_j, Y_j}\}_{j = 1}^m$ is independent. 
    
    \begin{proof}
        Induct on $m$. The base case of $m = 0$ holds by Corollary \ref{cor:elem-indep} and the inductive step works by applying the subtraction rule as follows. Let $S = \bigcap_{i \in I} S_{U_i}^{X_i} \cap \bigcap_{j=1}^{m} S_{V_j}^{Y_j, Z_j}$. Then: 
        \begin{align*} \mu\left(S \cap S_{V_{m+1}}^{Y_{m+1}, Z_{m+1}}\right) &= \mu\left(S \cap S_{V_{m+1}}^{Y_{m+1}}\right) - \mu\left(S \cap S_{V_{m+1}}^{Z_{m+1}}\right) = \mu(S) \left[\mu\left(S_{V_{m+1}}^{Y_{m+1}}\right) - \mu\left(S_{V_{m+1}}^{Z_{m+1}}\right)\right] \\ &= \mu\left(S\right) \mu\left(S_{V_{m+1}}^{Y_{m+1}, Z_{m+1}}\right).\end{align*}
    \end{proof}
\end{corollary}

Taken together, these properties reduce density computation to finding a partition of the set in question into intersections of simple sets supported on disjoint subschemes of $\Pp^n$. In this way, our problem is transformed into one of enumerating possible configurations of first-order local constraints that will guarantee that a curve is transverse-free. In the next section we put this basic idea into practice, obtaining simple results that are suggestive of the truth of the paper's title.

While it is unclear whether $\mu$ is a true probability measure defined on an appropriate $\sigma$-algebra, it shares enough similar properties that it is convenient to interpret it as one. It therefore makes sense to define a notion of \textbf{conditional density} in analogy with conditional probability. More precisely, we will denote \[\mu(A | B) \coloneqq \frac{\mu(A \cap B)}{\mu(B)}.\] The intuition provided by conditional probability will play an important role in Section \ref{sec:Sam-Sung’s-Note}. It will suggest a way to avoid cumbersome case-by-case calculations that would otherwise arise when decomposing $\Omega$ into simple sets. From the independence of simple sets we have a similar result for the conditional version of $\mu$. 

\begin{corollary}\label{cor:cond-indep}
    Suppose $\{U_i\}_{i=1}^m$ are mutually disjoint. If $B_i$, and $A_i \cap B_i$, are all simple sets supported at $U_i$, then we have that \[\mu\left(\bigcap_{i=1}^m A_i \bigg| \bigcap_{i=1}^m B_i\right) = \prod_{i=1}^m \mu(A_i | B_i) = \prod_{i=1}^m \mu\left(A_i \bigg| \bigcap_{i=1}^m B_i\right).\]

    \begin{proof}\begin{align*}&\mu\left(\bigcap_{i=1}^m A_i \bigg| \bigcap_{i=1}^m B_i\right) = \frac{\mu(\bigcap_{i=1}^m (A_i \cap B_i))}{\mu(\bigcap_{i=1}^m B_i)} = \prod_{i=1}^m \frac{\mu(A_i \cap B_i)}{\mu(B_i)} = \prod_{i=1}^m \mu(A_i | B_i) \\ &= \prod_{i=1}^m \frac{\mu(A_i \cap B_i)\mu(\bigcap_{j \neq i} B_j)}{\mu(B_i)\mu(\bigcap_{j \neq i} B_j)} = \prod_{i=1}^m \frac{\mu(A_i \cap \bigcap_{i=1}^m B_i)}{\mu(\bigcap_{i=1}^m B_i)} = \prod_{i=1}^m \mu\left(A_i \bigg| \bigcap_{i=1}^m B_i\right).\end{align*}\end{proof}
\end{corollary}

\subsection{Rational Density Values}\label{subsec: rationality}

We now characterize the subsets of $R_\homog$ whose densities can be computed using the techniques already introduced. It is at this point appropriate to make some remarks regarding subschemes of $\Pp^n$. Since $\Pp^n$ is Noetherian, we can freely interchange the order of a composition of a closed embedding and open embedding into $\Pp^n$. Now define for a subscheme $X$ and a closed subscheme $Y$, the open subscheme $X \setminus Y \coloneqq X - X \cap Y$ of $X$. It follows that every subscheme of $U \subseteq \Pp^n$ is quasiprojective and there exist closed subschemes of $Z_1, Z_2 \subseteq \Pp^n$ so that $U = Z_1 \setminus Z_2$. One may have already noticed that the definition of $S_U^X$ depends only on the topological subset underlying $U$. And indeed, every subscheme corresponds to a locally closed subset while every locally closed subset corresponds to a unique reduced subscheme. 

Let $\mathfrak{S}$ be the algebra over $R_\homog$ generated by simple sets, that is, the family of sets that can be formed by taking complements, finite unions, and finite intersections of simple sets. Suppose $\mathcal{P} \in \mathfrak{S}$ is generated by simple sets supported at subschemes that are all members of a finite collection $\mathcal{U}$. We will extend the notion of support to say that $\mathcal{P}$ is \textbf{supported} on every \textit{topological subset} $W$ which contains the underlying topological subset of every $U \in \mathcal{U}$. Due to the correspondence between subschemes and locally closed topological subsets, for our convenience we will subsequently make no special effort to distinguish between them and just refer to supports. We need the following lemma regarding the structure of members of $\mathfrak{S}$. 

\begin{lemma}\label{lem:structure}
    For every $\mathcal{P} \in \mathfrak{S}$ supported on $W$, there exists a finite collection, $\mathcal{U}$, of disjoint subschemes of $\Pp^n$ that are contained in $W$, a non-negative integer $m$, and non-negative integers $m_U$ for all $U \in \mathcal{U}$ so that $\mathcal{P}$ can be written as a finite disjoint union of finite intersections of simple sets, \[\mathcal{P} = \bigsqcup_{i=1}^m \left(\bigcap_{U \in \mathcal{U}} \bigcap_{j=1}^{m_U} S_{U, i, j}\right),\] where each $S_{U,i,j}$ is a simple set of the form $S_U^X\textrm{ or }S_U^{X, Y}$.

    \begin{proof}
        For any such $\mathcal{P}$ there is already a collection $\mathcal{U}$ of subschemes of $\Pp^n$ contained in $W$ so that we have a formula for $\mathcal{P}$ in terms of simple sets with lower entries in $\mathcal{U}$. Complements of simple sets are always a finite disjoint union of simple sets due to the identities $(S_U^X)^\complement = S_U^{\emptyset, X}$ and $(S_U^{X, Y})^\complement = S_U^{\emptyset, X} \sqcup S_U^Y$, so standard Boolean algebra manipulations will obtain an expression for $\mathcal{P}$ in the desired form with the same $\mathcal{U}$. The goal is to replace $\mathcal{U}$ with a disjoint collection by taking advantage of the subdividing identity $S_U^X = S_{U \cap Z}^X \cap S_{U \setminus Z}^X$ for $Z \subseteq \Pp^n$ closed. From it we obtain an analogous subdividing identity for complementary sets:
        \begin{align*}
            S_U^{X, Y} &= S_U^X \setminus (S_{U \cap Z}^Y \cap S_{U \setminus Z}^Y) = (S_U^X \setminus S^Y_{U \setminus Z}) \sqcup ((S_U^X \cap S_{U \setminus Z}^Y) \setminus S_{U \cap Z}^Y)) \\ &= (S_{U \cap Z}^{X, Y} \cap S_{U \setminus Z}^Y) \sqcup (S_{U \cap Z}^X \cap S_{U \setminus Z}^{X, Y}).
        \end{align*}
        Now suppose $U, U' \in \mathcal{U}$ are distinct but fail to be disjoint. Then there exists closed subschemes $Z_1, Z_2 \subseteq \Pp^n$ so that $U' = Z_1 \setminus Z_2$ and therefore $S_{U, i, j}$ can be written as the disjoint union of intersections of simple sets whose lower entries are in the collection $\{U \setminus Z_1, (U \cap Z_1) \setminus Z_2, U \cap Z_1 \cap Z_2\}$. These are clearly disjoint and importantly $(U \cap Z_1) \setminus Z_2 = U \cap U'$. Therefore, it is in fact possible to repeatedly apply subdividing identities until the lower entries form a collection whose members can be used to form partitions of each $U \in \mathcal{U}$. 
    \end{proof}
\end{lemma}

Say that $\mathcal{P} \subseteq R_{\homog}$ \textbf{splits} if we can write $\mathcal{P}$ in the form given by the previous lemma where $m_U = 1$ for all $U \in \mathcal{U}$. Evidently if $\mathcal{P}$ splits then it is also in $\mathfrak{S}$. This splitting property has two important consequences. First, if $\mathcal{P}$ splits then $\mu(\mathcal{P}) \in \Q$ and can be computed by applying the additivity of $\mu$, Corollary \ref{cor:simple-indep}, and Corollary \ref{cor:zeta}. Secondly, Corollary \ref{cor:simple-indep} can be upgraded as we have that subsets of $R_\homog$ that split and have disjoint supports are independent. So we now wish to compare $\mathfrak{S}$ with the set of its members that split, the obstruction here being the failure of intersections of the form $\bigcap_{j=1}^{m_U} S_{U,i,j}$ to simplify into a disjoint union of simple sets supported at $U$. The next lemma indicates when this cannot occur given $U$.

\begin{lemma}\label{lem:splitting}Let $U$ be a subscheme of $\Pp^n$ with $\dim{U} \in \{0, n-1, n\}$. Then there exists a finite partition $R_\homog$ into sets that split and are supported at $U$, $R_\homog = \bigsqcup S_i$, having the property that every simple set of the form $S_U^X\textrm{ or }S_U^{X,Y}$ has a partition obtained from a subset of the $\{S_i\}$. 
    \begin{proof} If $U = \emptyset$, then $R_\homog$ is a suitable partition of itself so assume otherwise. We have two cases: \begin{enumerate}
        \item $\dim{U} = 0$: Then $U$ is the disjoint union of a finite number of closed points. Since $S_P^{X, Y} = S_P^X \setminus S_P^Y$, we may reduce to constructing an appropriate partition of $R_\homog$ with respect to just the elementary sets $S_P^X$, and because $S_{U_1 \coprod U_2}^X = S_{U_1}^X \cap S_{U_2}^X$, we may further reduce to the case where $U$ is a single closed point $P$.

        We now refer to the material in \cite{poonen} in an aim to generalize the  picture provided in Lemma $\ref{lem:numbers}$. Suppose $X$ is a smooth subscheme of $\Pp^n$ containing $P$. Then we may let $\mathfrak{m}_X$ be the ideal sheaf of $P$ on $X$ and $Y_X$ be the closed subscheme of $X$ corresponding to the ideal sheaf $\mathfrak{m}_X^2 \subseteq \Orb_X$. Now consider the special case $Z = Y_{\Pp^n}$. There will be some coordinate $x_i$ that is invertible on $Z$ and as such, for $f \in R_d$ we will let $f\vert_Z \in H^0(Z, \Orb_Z)$ be the restriction of $x_i^{-d}f$ to $Z$. Considering the inclusion structure on smooth subschemes of $\Pp^n$ we see that we have a chain of maps \[H^0(Z, \Orb_Z) \xrightarrow{\phi_X} H^0(Y_X, \Orb_{Y_X}) \to H^0(Y_P, \Orb_{Y_P}).\] In particular, for $f \in R_\homog$, $H_f \cap X$ is smooth of dimension $\dim{X} - 1$ when $\phi_X(f\vert_Z) \in H^0(P_X, \Orb_{Y_X})$ is non-zero. 
        
        We additionally know that the $H^0(P_X, \Orb_{Y_X})$ are vector spaces over the residue field $\kappa(P)$ of dimension $\dim X + 1$ and the maps $\phi_X$ are  $\kappa(P)$-linear. Correspondingly, we can construct an equivalence relation on $H^0(Z, \Orb_Z)$ where two elements $g_1, g_2$ are considered equivalent whenever there exists $u \in \kappa(P)^\times$ so that $g_1 = ug_2$. Denote the equivalence class of $g \in H^0(Z, \Orb_Z)$ under this relation by $[g]$. Then  for any $X$, whether $f \in S_U^X$ is completely determined by $[f\vert_Z]$.
        
        So now let $\mathcal{P}_{[g]}$ be the subset of $R_\homog$ for which $f\vert_Z \in [g]$. Due to the conditions for smooth intersection with $\Pp^n$ and $P$, we see that $R_\homog$ can be partitioned as \[R_\homog = S_{P}^{\emptyset, \Pp^n} \sqcup \left(\bigsqcup_{[g] \subseteq \phi_P^{-1}(\{0\}) \setminus [0]} \mathcal{P}_{[g]}\right) \sqcup S_P^P.\] In particular if $f \in S_{P}^{\emptyset, \Pp^n}$ then $f \notin S_P^X$ for any $X$, and if $f \in S_{P}^{P}$, then $f \in S_P^X$ for any $X$. So it suffices to show that these sets $\mathcal{P}_{[g]}$ are in fact simple sets supported at $P$.

        But now we may apply Theorem 1.2 of \cite{poonen}, to show that for every equivalence class $[g] \neq [0]$, there exists $\tilde{g} \in R_\homog$ so that $\tilde{g}\vert_Z \in [g]$ and $H_{\tilde{g}}$ is a smooth hypersurface. We know that $[g]$ is in the kernel of $\phi_{H_{\tilde{g}}}$ and since $\dim_{\kappa(P)} H^0(Z, \Orb_Z) - \dim_{\kappa(P)} H^0(P_X, \Orb_{Y_X}) = 1$, it follows that actually $\ker(\phi_{H_{\tilde{g}}}) = [g] \sqcup [0]$. Equivalently, $\mathcal{P}_{[g]} \sqcup \mathcal{P}_{[0]} = S_P^{\emptyset, H_{\tilde{g}}}$ and it follows that $\mathcal{P}_{[g]} = S_P^{\Pp^n, H_{\tilde{g}}}$ as desired. Therefore, there  exists a finite collection of hypersurfaces containing $P$, $\mathcal{H}$, such that we have an appropriate partition $R_\homog = S_{P}^{P} \sqcup \left(\bigsqcup_{H \in \mathcal{H}} S_{P}^{\Pp^n, H}\right) \sqcup S_{P}^{\emptyset, \Pp^n}$.
        
        \item $\dim{U} \in \{n-1, n\}$: Consider $R_\homog = S_U^U \sqcup S_U^{\Pp^n, U} \sqcup S_U^{\emptyset, \Pp^n}$. It is sufficient to check whether the requisite partitions exist just for the elementary sets $S_U^X$. The case $X = \emptyset$ is clear so we may otherwise assume $U \subseteq X \subseteq \Pp^n$. Then either $\dim X = \dim U$ or $\dim X = \dim{\Pp^n}$. In the former case, since we require $X$ smooth, $U \subseteq X$ is open and therefore $S_U^X = S_U^U$. In the latter case, $X \subseteq \Pp^n$ open so $S_U^X = S_U^{\Pp^n} = S_U^U \sqcup S_U^{\Pp^n, U}$.
    \end{enumerate}
    \end{proof}
\end{lemma}

\begin{theorem}\label{thm:inQ}
If $n \leq 2$, then for all $\mathcal{P} \in \mathfrak{S}$, $\mathcal{P}$ splits, and thus, $\mu(\mathcal{P}) \in \Q$. 
    \begin{proof}Since $\mathcal{P} \in \mathfrak{S}$, we can write $\mathcal{P}$ in the form given by \ref{lem:structure} where we do not necessarily have $m_U = 1$ for all $U \in \mathcal{U}$. Since $n \leq 2$, we know that $\dim{U} \in \{0, n-1, n\}$ for any $U \subseteq \Pp^n$, and therefore the previous lemma guarantees that $\bigcap_{j=1}^{m_U} S_{U,i,j}$ has a finite partition into subsets of $R_\homog$ that split and are supported at $U$. Then distributing the $\cap$'s over the $\sqcup$'s gives the result. \end{proof}
\end{theorem}

\begin{remark}
     If $n \geq 3$, the set of hypersurfaces that intersect two distinct hyperplanes transversely is captured by $S_{\pi_1}^{\pi_1} \cap S_{\pi_2}^{\pi_2} = S_{\pi_1 \setminus \pi_2}^{\pi_1} \cap S_{\pi_2 \setminus \pi_1}^{\pi_2} \cap \left(S_{\pi_1 \cap \pi_2}^{\pi_1} \cap S_{\pi_1 \cap \pi_2}^{\pi_2}\right)$. Therefore if the density of this set is computable, we should be able to compute $\mu(S_{\pi_1 \cap \pi_2}^{\pi_1} \cap S_{\pi_1 \cap \pi_2}^{\pi_2})$. However, assuming the arguments in \cite{poonen} generalize appropriately, we obtain \[\mu(S_{\pi_1 \cap \pi_2}^{\pi_1} \cap S_{\pi_1 \cap \pi_2}^{\pi_2}) = \prod_{\textrm{closed }P \in \pi_1 \cap \pi_2} (1 - 2q^{n \deg{P}}),\] which is an infinite product that numerical estimates do not suggest has a rational value.
\end{remark}

We end the section by returning to our comment that subsets of $R_\homog$ that split and have disjoint supports are independent. We can upgrade Corollary \ref{cor:cond-indep} immediately: 

\begin{corollary}
    Suppose $\{U_i\}_{i=1}^m$ are mutually disjoint. If $B_i$ splits and is supported on $U_i$, while $A_i \cap B_i$ splits and is supported on $U_i$ for every $i$, then \[\mu\left(\bigcap_{i=1}^m A_i \bigg| \bigcap_{i=1}^m B_i\right) = \prod_{i=1}^m \mu(A_i | B_i) = \prod_{i=1}^m \mu\left(A_i \bigg| \bigcap_{i=1}^m B_i\right).\]
\end{corollary}

For the sake of stating and proving Theorem \ref{thm:Sam-Sung’s-Note} later on, we are interested in the situation where the $B_i$ are identically some $B \subseteq R_\homog$, and we rely only on the disjointness of the supports of the $A_i$ to provide independence. This situation may be dealt with provided we introduce one additional piece of terminology.

We say that $\mathcal{P} \subseteq R_\homog$ \textbf{factors through} a finite disjoint collection of supports, $\mathcal{U}$, if $\mathcal{P} = \bigcap_{U \in \mathcal{U}} S_U$ where each $S_U$ splits and is supported on $U$. Notice that $\mathcal{P}$ also factors through any partition of $\bigsqcup_{U \in \mathcal{U}} U$ coarser than $\mathcal{U}$, while if the $S_U$ are all elementary sets, then $\mathcal{P}$ also factors through any partition of $\bigsqcup_{U \in \mathcal{U}} U$ finer than $\mathcal{U}$. As our subsequent work will only deal with the case $n = 2$, we take advantage of Theorem \ref{thm:inQ} and obtain the further corollary:

\begin{corollary}\label{cor:buffed}
    Assume $n \leq 2$ and suppose $\{U_i\}_{i=1}^m$ are mutually disjoint with $A_i \in \mathfrak{S}$ supported on $U_i$. If $B$ factors through $\mathcal{U}$ that is finer than the collection $\{U_i\}_{i=1}^m$ in the sense that for every $U \in \mathcal{U}$, there exists at most one $i$ such that $U \cap U_i \neq \emptyset$, then
    \[\mu\left(\bigcap_{i=1}^m A_i \bigg| B\right) = \prod_{i=1}^m\mu(A_i | B).\]

    \begin{proof}
        We obtain a partition $\{U'_i\}_{i=1}^{m+1}$ of $\bigcup_{U \in \mathcal{U}} U$ coarser than $\mathcal{U}$ by letting $U'_i$ be the union of all $U \in \mathcal{U}$ such that $U \cap U_i \neq \emptyset$ except for when $i = m+1$, in which case we let $U'_{m+1}$ be the union of all $U \in \mathcal{U}$ such that $U \cap U_i = \emptyset$ for all $i$. Additionally set $A_{m+1} = S_{U'_{m+1}}^\emptyset = R_\homog$. We have $U_i \subseteq U'_i$, so in fact $A_i$ is always supported on $U'_i$. We can then apply the previous corollary to obtain the desired result using these $\{A_i\}_{i=1}^{m+1}$ supported on the disjoint $\{U'_i\}_{i=1}^{m+1}$ with the $\{B_i\}_{i=1}^{m+1}$ constructed from the fact that $B$ factors through $\{U'_i\}_{i=1}^{m+1}$.
    \end{proof}
\end{corollary}

\section{Trivially Transverse-Free Curves}\label{sec:TTFC}

In contrast to the situation explored in \cite{freidin-asgarli}, our curves are allowed to have singularities and distinct lines may fail to have a transverse intersection at the same point on the curve. A curve that is singular at a rational point will fail to transversely intersect every line passing through it. Therefore, any blocking set corresponds to a class of transverse-free curves with singularities at each point of the blocking set.

As noted earlier, the simplest blocking sets are the trivial blocking sets containing the $q+1$ rational points of a line. So for a line $\ell$, let \[\Theta_{\ell} \coloneqq \{f \in \R_\homog|C_f \text{ has singularities at every rational point of } \ell\}.\] We may then define trivially transverse-free curves as those curves whose set of rational singular points is a trivial blocking set. The corresponding set of polynomials is \[\Theta \coloneqq \bigcup_{\ell \subset \Pp^2} \Theta_{\ell} = \{f \in R_\homog | C_f \text{ is trivially transverse-free}\}.\] The inclusion-exclusion principle gives a tight bound on $\mu(\Theta)$, and hence a lower bound for $\mu(\Omega)$:

\begin{theorem}\label{thm:lubber}
     \[(q^2 + q + 1)q^{-3(q+1)}(1 - q^{-3q+2}) \leq \mu(\Theta) \leq (q^2+q+1)q^{-3(q+1)}.\] 
\end{theorem}

\begin{proof}
Truncating the expansion given by the inclusion-exclusion principle after the first term will give an upper bound while truncating after the second term will give a lower bound. So iterating over lines in $\Pp^2$, we have \[\sum_{\ell} \mu(\Theta_{\ell_i}) - \sum_{\ell_1 \neq \ell_2} \mu(\Theta_{\ell_1} \cap \Theta_{\ell_2}) \leq \mu(\Theta) \leq \sum_{\ell} \mu(\Theta_{\ell}).\] Every line contains exactly $q+1$ rational points, and every two distinct lines intersect at exactly $1$ rational point so we may apply Corollary \ref{cor:simple-indep} and Lemma \ref{lem:numbers} to obtain \begin{align*}\mu\left(\Theta_{\ell}\right) &= \mu\left(\bigcap_{P \in \ell(\Fq)} S_P^{\emptyset, \Pp^2}\right) = \left(q^{-3}\right)^{|\ell(\Fq)|} = q^{-3(q+1)}, \\ \mu\left(\Theta_{\ell_1} \cap \Theta_{\ell_2}\right) &= \mu\left(\bigcap_{P \in \ell_1(\Fq) \cup \ell_2(\Fq)} S_P^{\emptyset, \Pp^2}\right) = \left(q^{-3}\right)^{|\ell_1(\Fq) \cup \ell_2(\Fq)|} = q^{-3(2q+1)}.\end{align*} Substituting this into the summation obtains \[(q^2+q+1)q^{-3(q+1)} - \binom{q^2+q+1}{2} q^{-3(2q+1)} \leq \mu(\Theta) \leq (q^2+q+1)q^{-3(q+1)},\] and we get the lower bound claimed by using $q \leq q^2$ and factoring.
\end{proof}

In the other direction, we might ask how close this lower bound is to an upper bound for $\mu(\Omega)$. We notice the heuristic that while specifying rational singular points `costs more' than merely specifying that it lies on the curve (compare $q^{-3}$ with $q^{-1}$ in Lemma \ref{lem:numbers}), it is vastly more `cost-effective', accounting for $q+1$ non-transverse lines, while the alternative only accounts for $1$. Such a rough `efficiency calculation' therefore suggests that for sufficiently large $q$, those curves with a blocking set of singularities should account for the bulk of the density of transverse-free curves. 

The example of trivially transverse-free curves suggests that focusing on just one of the lines in the plane may be enough to provide a simple upper bound. Then every other line intersects it somewhere and if $C_f$ is not trivially transverse-free, we should be able to leverage the existence of rational points on this line that are not singularities of $C_f$ to our advantage. We will be interested in the density provided by the following lemma.

\begin{lemma}\label{lem:resurrected}
    \[\mu(\{f \in R_\homog : \text{ $C_f$ has a non-transverse intersection with $\ell$ at $P' \neq P$}\}) = q^{-1}.\]

    \begin{proof}
        This set is the same as $S_{\ell \setminus P}^{\emptyset, \ell}$. We apply Corollary \ref{cor:simple-indep} in reverse to get $\mu(S_{\ell \setminus P}^{\emptyset, \ell}) = \mu(S_{\ell \setminus P}^\emptyset) - \mu(S_{\ell \setminus P}^\ell) = 1 - \frac{\mu(S_{\ell \setminus P}^\ell \cap S_P^\ell)}{\mu(S_P^\ell)} = 1 - \frac{\mu(S_{\ell}^\ell)}{\mu(S_P^\ell)}$  and so Corollary \ref{cor:zeta} finishes.
    \end{proof}
\end{lemma}

It enables us to obtain the following rudimentary upper bound.

\begin{theorem}\label{thm:old-ubber}
    \[\mu (\Omega) \leq 2q^{-(q+1)}.\]  
\end{theorem}

\begin{proof}
    Our goal is to cover $\Omega$ by a family of sets whose densities we can compute. To do this we notice that every transverse-free curve must fail to transversely intersect at least the $q + 1$ lines passing through some rational point, $P$. These lines only intersect each other at $P$ so with a little work we will be in a position to apply independence. Pick a line, $\ell_0$, and label its rational points as $P_0, P_1, \dots, P_q$. For any curve, $C_f$, exactly one of the following three statements must be true: 
    
    \begin{enumerate}
        \item $C_f$ is not singular at $P_0$.
        \item There exists $1\leq k \leq q$ such that $C_f$ is singular at $P_0, P_1, \dots, P_{k-1}$ but not singular at $P_k \in \ell_0$.
        \item $C_f$ is singular at every rational point of $\ell_0$.
    \end{enumerate} 

    For the first case, let the set of lines passing through $P_0$ be $L$. We have that if $f \in \Omega$, then either $C_f$ does not contain $P_0$ and $C_f$ is has a non-transverse intersection with each line in $L$ at a point that is not $P_0$, or $C_f$ does not intersect exactly one of the lines in $L$, transversely at $P_0$, and has a non-transverse intersection with each of the remaining lines in $L$ at a point that is not $P_0$. This organization of constraints on $C_f$ allows all of Lemma \ref{lem:numbers}, Lemma \ref{lem:resurrected}, and Corollary \ref{cor:simple-indep} to be applicable and we have:  
    
    \[\mu(\{f \in \Omega \colon C_f \textrm{ is not singular at } P_0\}) \leq \mu\left(\left(S_{P_0}^{P_0} \cap \bigcap_{\ell \in L} S_{\ell \setminus P_0}^{\emptyset, \ell}\right) \sqcup \bigsqcup_{\ell' \ni P_0} \left(S_{P_0}^{\Pp^2, \ell} \cap \bigcap_{\ell \in L \setminus \{\ell'\}} S_{\ell \setminus P_0}^{\emptyset, \ell}\right)\right).\]
    Further,
    \begin{align*}
    &= \mu(S_{P_0}^{P_0})\mu(S_{\ell \setminus P_0}^{\emptyset, \ell})^{|L|} + (q+1) \cdot \mu(S_{P_0}^{\Pp^2, \ell})\mu(S_{\ell \setminus P_0}^{\emptyset, \ell})^{|L \setminus \{\ell'\}|} \allowdisplaybreaks \\&\leq \mu(S_{P_0}^{P_0})\mu(S_{\ell \setminus P_0}^{\emptyset, \ell})^{|L|} + (q+1) \cdot \mu(S_{P_0}^{\Pp^2, \ell})\mu(S_{\ell \setminus P_0}^{\emptyset, \ell})^{|L| - 1} = (2 - q^{-1} - q^{-2})(q^{-1})^{|L|}.\end{align*}
    
    Similarly, for the second case, let the set of lines passing through $P_k$ that aren't $\ell_0$ be $L_k$. Again if $f \in \Omega$, then every line in $L_k$ is non-transverse, and the breakdown of possibilities in the previous case still holds with the added stipulation that $C_f$ is also singular at $P_0, P_1, \dots, P_{k-1}$. Because none of the lines in $L_k$ pass through $P_0$, independence still holds with $S_{P_i}^{\emptyset, \Pp^2}$ for $0 \leq i \leq k-1$ and so we obtain analogously \begin{align*}&\mu(\{f \in \Omega \colon C_f \textrm{ is singular at } P_0, P_1, \dots, P_{k-1} \textrm{ and not singular at } P_k\}) \\&\leq (2 - q^{-1} - q^{-2})(q^{-1})^{|L_k|}\prod_{i=0}^{k-1}\mu(S_{P_i}^{\emptyset, \Pp^2}) = (2 - q^{-1} - q^{-2})(q^{-1})^{|L| - 1}(q^{-3})^k.\end{align*}
    
    The third case we know already from the proof of the previous theorem. It is $\mu(\Theta_{\ell_0}) = (q^{-3})^{q+1}$. Now $\mu(\Omega)$ is bounded above by the sum of the upper bounds computed for each possibility outlined: \[[1 + (q^{-3} + (q^{-3})^2 + \dots + (q^{-3})^q) \cdot (q^{-1})^{-1}](2 - q^{-1} - q^{-2})(q^{-1})^{|L|}  + (q^{-3})^{q+1}.\]
    
    Now $1 + (q^{-3} + (q^{-3})^2 + \dots + (q^{-3})^q) \cdot (q^{-1})^{-1} = 1 + q^{-2} + q^{-5} + \dots + q^{-(3q-1)} \leq 1 + q^{-2} + q^{-4}$, so after using $|L| = q + 1$ we find that the claimed bound results from a slight weakening: \[[1 + q^{-2} + q^{-4}](2 - q^{-1} - q^{-2})(q^{-1})^{q+1} + q^{-3(q+1)} \leq (2 - q^{-2} + q^{-2q-2})q^{-(q+1)} \leq 2q^{-(q+1)}.\]
\end{proof}

The proof above suggests the following strategy. First, prescribe the $\Fq$-rational singularities of the curve in advance. Then quantify how the lines avoiding these singularities further constrain the density of possible transverse-free curves in each case. The main difficulty here is that these lines are not independent, since the lines inevitably intersect. The next section introduces a sub-independence result which overcomes this difficulty, paving the way for a generalized version of this proof that attains the much better explicit bound in Theorem \ref{thm:upper-bound}.

\section{
Sub-Independence of Non-Transverse Lines}\label{sec:Sam-Sung’s-Note}

Using the notation we developed, $\Omega$ can be cleanly expressed as $\bigcap_\ell S_\ell^{\emptyset,\ell}$. However these simple sets are not independent because any two lines intersect at a rational point and direct computation is otherwise intractable. On the other hand, the product of the individual densities is still a natural candidate for an upper bound on density. This fails if tried directly as one can check that the lower bound we have already obtained is incompatible with this strategy. The trick here is to massage the situation into one in which this works. We now present a lemma to encapsulate this idea.

In particular, we provide an upper bound that can be interpreted as stating that in the absence of singularities, the events of being not transverse to different lines are negatively related. Indeed, a curve cannot have a non-transverse intersection with each of two lines at the point common to both lines without having a singularity there. The main difficulty in the proof is extending this idea to subsets of $R_\homog$ that have more complicated descriptions in terms of simple sets. We need to ensure that densities remain tractably computable throughout. To handle this, we use the notions defined in Section \ref{subsec: rationality}.

Recall that if $\mathcal{P} \subseteq R_\homog$ can be written in terms of complements, finite unions, and finite intersections of simple sets supported at subschemes $U_i$, then $\mathcal{P}$ is supported at $U$ whenever $U$ contains all of the $U_i$. As a consequence of Theorem \ref{thm:inQ}, we can say that $\mathcal{P} \subseteq R_\homog$ factors through a disjoint collection of supports, $\mathcal{U}$, whenever $\mathcal{P} = \bigcap_{U \in \mathcal{U}} S_U$ for $S_U \subseteq R_\homog$ that is supported on $U$ in the manner just exhibited. This notion of factoring is the one required to apply Corollary \ref{cor:buffed} guaranteeing the independence of conditional densities.

\begin{theorem}\label{thm:Sam-Sung’s-Note}
    Let $L$ be a finite set of lines and suppose that $\Psi \subset R_\homog$ is such that for all $f \in \Psi$, $C_f$ is smooth at any intersection point of two distinct lines in $L$. If we also have that $\Psi$ factors through $\mathcal{U}$ for which every $U \in \mathcal{U}$ satisfies either $U \cap \ell = \emptyset$ or $U \subseteq \ell$ for every $\ell \in L$, then: \[\mu\left(\bigcap_{\ell \in L} S_\ell^{\emptyset,\ell} \bigg|\Psi\right) \leq \prod_{\ell \in L} \mu(S_\ell^{\emptyset,\ell}|\Psi).\]
    \begin{proof}
        The theorem reduces to showing that for any $L' \subseteq L$, and $\ell_0 \in L \setminus L'$, \[\mu\left(\bigcap_{\ell \in L' \cup \{\ell_0\}}S_\ell^{\emptyset,\ell}\bigg|\Psi\right) \leq \mu\left(\bigcap_{\ell \in L'}S_\ell^{\emptyset,\ell}\bigg|\Psi\right)\mu(S_\ell^{\emptyset,\ell}|\Psi),\] because the result would then follow by inducting on the cardinality, $|L|$. The key here is choosing an appropriate partitioning of these conditional densities so that the claim follows on each of the pieces via using the independence result provided in Corollary \ref{cor:buffed}. For any $M \subseteq L'$, let \[T_M \coloneqq \{f \in R_\homog | \text{$C_f$ fails to transversely intersect $\ell \in L'$ at a point not on $\ell_0$ iff $\ell \notin M$}\}.\]  Symbolically, $T_M = \left(\bigcap_{\ell \in M}S_{\ell \setminus \ell_0}^\ell\right) \cap \left(\bigcap_{\ell \in L' \setminus M}S_{\ell \setminus \ell_0}^{\emptyset, \ell}\right)$ and by construction, the $T_M$ give a partition of $R_\homog$. We note that $T_M$ is supported on $\Pp^2 \setminus \ell_0$. 
        
        Now if a curve $C_f \in  T_M  \cap \left(\bigcap_{\ell \in L'}S_\ell^{\emptyset,\ell}\right)$, then it must be true that for each $\ell \in M$, $C_f$ does not intersect $\ell$ transversely at the point $\ell \cap \ell_0$ since the required failure point can't be elsewhere. As such, $T_M  \cap \left(\bigcap_{\ell \in L'}S_\ell^{\emptyset,\ell}\right) = T_M \cap \left(\bigcap_{\ell \in M}S_{\ell \cap \ell_0}^{\emptyset, \ell}\right)$. By the smoothness constraint on $\Psi$, it follows that a curve $C_f \in  T_M  \cap \left(\bigcap_{\ell \in L'}S_\ell^{\emptyset,\ell}\right) \cap \Psi$ is in $S_{\ell_0}^{\emptyset,\ell_0}$ if and only if it is in $S_{\ell_0 \setminus \bigcup_{\ell \in M} \ell}^{\emptyset,\ell_0}$. 
        
        So now notice that $S_{\ell_0 \setminus \bigcup_{\ell \in M} \ell}^{\emptyset,\ell_0}$ is supported on $\ell_0 \setminus \bigcup_{\ell \in M} \ell$ while $S_{\ell \cap \ell_0}^{\emptyset, \ell} \cap T_M$ is supported on $\Pp^2 \setminus \left(\ell_0 \setminus \bigcup_{\ell \in M} \ell\right)$. These two supports form a partition of $\Pp^2$ and our requirement on $\Psi$ guarantees that it factors through some $\mathcal{U}$ whose members must be strictly contained in one of these. It follows that Corollary \ref{cor:buffed} may be applied using the resulting set equalities to obtain
        \begin{align*}&\mu\left(\left(\bigcap_{\ell \in L' \cup \{\ell_0\}}S_\ell^{\emptyset,\ell}\right) \cap T_M \Bigg| \Psi\right) = \mu\left(S_{\ell_0 \setminus \bigcup_{\ell \in M} \ell}^{\emptyset,\ell_0} \cap \left(\bigcap_{\ell \in M}S_{\ell \cap \ell_0}^{\emptyset, \ell}\right) \cap T_M \bigg| \Psi\right) \\ &= \mu\left(\left(\bigcap_{\ell \in M}S_{\ell \cap \ell_0}^{\emptyset, \ell}\right) \cap T_M \bigg| \Psi\right)\mu\left(S_{\ell_0 \setminus \bigcup_{\ell \in M} \ell}^{\emptyset,\ell_0} \bigg| \Psi\right)  \leq \mu\left(\left(\bigcap_{\ell \in L'}S_\ell^{\emptyset,\ell}\right) \cap T_M \bigg| \Psi\right)\mu\left(S_{\ell_0}^{\emptyset,\ell_0} \bigg| \Psi\right).\end{align*} Now by summing over all $M \subseteq L'$ we obtain the desired reduction.
    \end{proof}
\end{theorem}

The upshot of this presentation is its versatility; we will now give three separate applications of it. For our purposes, $\Psi$ will always factors as the intersection of elementary sets as well as the intersection of complementary sets supported at a closed points, so the requirement that $\Psi$ factors through a sufficiently fine collection of supports is immediate. This is because we can just keep applying the subdividing identity $S_U^X = S_{U \cap \ell}^X \cap S_{U \setminus \ell}^X$ for elementary sets as necessary.

For our main question we can apply Theorem \ref{thm:Sam-Sung’s-Note} as follows. For $S \subseteq \Pp^2(\Fq)$, let $M_S$ be the set of lines that pass through no points of $S$ and let \[\Sing_S = \{f \in R_\homog | \text{The set of rational singular points of $C_f$ is $S$}\}.\]

\begin{lemma}\label{lem:diamond-prob}
    \[\mu(\Omega \cap \Sing_S) \leq (q^{-3})^{|S|}(1-q^{-3})^{q^2+q+1-|S|}\left(1 - \frac{(1-q^{-1})(1-q^{-2})}{(1-q^{-3})^{q+1}}\right)^{|M_S|} \leq \left(\frac{1}{q}\right)^{3|S|+|M_S|}.\]

    \begin{proof}
        If $C_f \in \Sing_S$, then any line passing through a point $P \in S$ will not be transverse to $C_f$. Therefore, we additionally have $C_f \in \Omega$ if and only if every $\ell \in M_S$ is not transverse to $C_f$. This establishes \[\mu(\Omega|\Sing_S) = \mu\left(\bigcap_{\ell \in M_S}S_\ell^{\emptyset,\ell}\bigg|\Sing_S\right).\] Now notice that $\Sing_S = \left(\bigcap_{P \in S} S_P^{\emptyset, \Pp^2}
        \right) \cap \left(\bigcap_{P \in \Pp^2(\Fq) \setminus S} S_P^{\Pp^2}\right)$, implying not only that $\mu(\Sing_S) = (q^{-3})^{|S|}(1-q^{-3})^{q^2+q+1-|S|}$, but also that any $C_f \in \Sing_S$ is smooth at any point where distinct lines in $M_S$ may intersect. Therefore, we can apply Theorem \ref{thm:Sam-Sung’s-Note} directly. It then suffices to compute $\mu(S_\ell^{\emptyset,\ell}|\Sing_S)$ for $\ell \in M_S$ to verify the first inequality, but we have that $1 - \mu(S_\ell^{\emptyset,\ell}|\Sing_S) = \mu(S_\ell^\ell|\Sing_S) = \mu(S_\ell^\ell|\bigcap_{P \in \ell(\Fq)} S_P^{\Pp^2}) = \mu(S_\ell^\ell)\mu(S_P^{\Pp^2})^{-(q+1)}$, which agrees. The last inequality follows since $1-q^{-3} \leq 1$ in the second factor and $1 - q^{-2} \geq (1-q^{-3})^{q+1}$ in the third factor.
    \end{proof}
\end{lemma}

Although somewhat unnecessary given that we can compute the next density exactly (see Theorem \ref{thm:old-ubber}), we can apply the theorem again except replacing $\Sing_S$ with \[\Sing'_{S, P} = \{f \in R_\homog | \text{$C_f$ is singular at every point in $S$, but not at $P$}\},\] where $P$ lies on every line in $M_S$. Explicitly, $\Sing'_{S, P} = S_P^{\Pp^2} \cap \left(\bigcap_{P \in S} S_P^{\emptyset, \Pp^2}\right)$, and analogously:

\begin{lemma}\label{lem:club-prob}
    \[\mu(\Omega \cap \Sing'_{S, P}) \leq (q^{-3})^{|S|}(1-q^{-3})\left(1 - \frac{(1-q^{-1})(1-q^{-2})}{(1-q^{-3})}\right)^{|M_S|} < \left(\frac{1}{q}\right)^{3|S|}\left(\frac{q+1}{q^2}\right)^{|M_S|}.\]
\end{lemma}

We save the easiest application of Theorem \ref{thm:Sam-Sung’s-Note} for last. Let the subset of $R_\homog$ corresponding to smooth \tf\, curves be $\Upsilon$. Then improving on \cite{freidin-asgarli}, we have:

\begin{theorem}\label{thm:smBounds}
     \[\mu(\Upsilon) \leq (1-q^{-1})(1-q^{-2})(1-q^{-3})\left(1-\frac{1-q^{-1}}{1-q^{-3}}\right)^{q^2+q+1}.\]

     \begin{proof}
         Straightforward. We have that $\Upsilon = S_{\Pp^2}^{\Pp^2} \cap \left(\bigcap_\ell S_\ell^{\emptyset,\ell}\right)$ and $S_{\Pp^2}^{\Pp^2}$ satisfies the necessary conditions on $\Psi$. Now $\mu(S_{\Pp^2}^{\Pp^2}) = (1-q^{-1})(1-q^{-2})(1-q^{-3})$ and $\mu(S_\ell^{\emptyset,\ell} | S_{\Pp^2}^{\Pp^2}) = 1 - \mu(S_\ell^\ell | S_{\Pp^2}^{\Pp^2}) = 1 - \mu(S_\ell^\ell | S_{\ell}^{\Pp^2}) = 1 - \mu(S_\ell^\ell)\mu(S_{\ell}^{\Pp^2})^{-1} = 1 - (1-q^{-1})(1-q^{-3})^{-1}$, so the claim follows by Theorem \ref{thm:Sam-Sung’s-Note}.
     \end{proof}
\end{theorem}
\section{ 
Enumerative Properties of Blocking Sets}\label{sec:block-time}
For more information on blocking sets the reader is encouraged to see the survey paper \cite{blocking-survey} on the spectrum problem. From a combinatorial perspective the projective plane is $\PG(2, q)$, the incidence structure consisting of the points and lines of $\Pp^2$ defined over $\Fq$. A blocking set is just a point set incident to all of the lines. In the ensuing discussion, let $k$ denote the size of a blocking set. For our purposes, the small minimal blocking sets are those containing less than $2q$ points (compare to the $\frac{3}{2}(q+1)$ standard in the literature). Then when $k$ is relatively small we have the following known results:  

From \cite{baer-equality} we have that minimal blocking sets of size at most $q + q^{\frac{1}{2}} + 1$ are either the $q + 1$ points of a line or the $q + q^{\frac{1}{2}} + 1$ points of a Baer subplane, the number of which is given in \cite{baer-number} as $(q^3 + q^{\frac{3}{2}})(q + 1)$. Stronger results may be obtained using the structural theorems due to Sziklai in \cite{best-bounds} which show that there exists an associated number called the  \textbf{exponent}, $e$, of a minimal blocking set that must divide $k$. Specifically we have that when $e = k$, the sets correspond to lines, when $e = k/2$, all are of size exactly $q + q^{\frac{1}{2}} + 1$ and  correspond to Baer subplanes, while those of exponent $e \leq k/3$  must have size at least $q + q^\frac{2}{3} + 1$. Recall that by $B_k$ we denote the number of minimal blocking sets of size $k$. To summarize: 

\begin{lemma}[\cite{baer-equality}, \cite{baer-number}, \cite{best-bounds}]\label{lem:tiny-B_k}
For integers $k < q + q^\frac{2}{3} + 1$, $B_k$ is known exactly. Namely:
    \[B_k = \begin{cases} 0 & \text{if } k \neq q+1, q + q^\frac{1}{2} + 1 \\
    q^2 + q + 1 & \text{if } k = q+1 \\
    (q^3 + q^\frac{3}{2})(q+1) & \text{if } k = q + q^\frac{1}{2} + 1 \end{cases}.\]
\end{lemma}

These quantities explain the observed estimates on $\mu(\Omega)$.  They are not however sufficient to determine $\mu(\Omega)$ to any reasonable accuracy. It turns out that it is necessary to have some control on $B_k$ up to somewhere just beyond the $k < \frac{3}{2}(q + 1)$ region for first order estimates, while the region $k < 2q$ is amenable to some simple analysis.

In that direction, Theorem \ref{theorem:min-blset-bound} obtains a very weak upper bound on $B_k$ for $k < 2q$ which to our knowledge is the first example of explicit bounds on $B_k$ for $\frac{3}{2}(q+1) < k < 2q$, and the only example for not too small $k$ that does not rely on any assumptions on the prime factorization of $q$. Its proof relies on exploiting a foundational result of Jamison, Brouwer, and Schrijver. Let $\mathrm{AG}(2, q)$ be the incidence structure formed by the points and lines of $\A^2 = \Pp^2 \setminus \ell$. Then: 
\begin{lemma}[Theorem 1 of \cite{affine-blocking-1}, Result 1.1 of \cite{affine-blocking-2}]\label{lem:affine-blset}
    A blocking set of $\AG(2,q)$ contains at least $2q-1$ points.
\end{lemma}

The original proofs of this fact exploit the algebraic structure unique to $\PG(2, q)$ as compared to other projective planes of order $q$, freeing the subsequent arguments in this paper to appeal to basic combinatorial structure instead. The importance of this result to our work cannot be overstated. In addition to providing upper bounds on $B_k$, it also accounts for the existence of non-trivial lower bounds on $|M_S|$, the number of lines not intersecting a set of points $S$. This quantity was introduced in Section \ref{sec:Sam-Sung’s-Note}, and these bounds explain why the contribution of non-blocking sets is negligible our estimates. The literature already provides the bounds we make use of, albeit in different forms. We reproduce their proofs here both to highlight their simplicity and clarify the exposition in the original sources.

\begin{lemma}[Corollary in \cite{small-blocking}]\label{lem:critical}
    If $S$ is a blocking set of $\PG(2, q)$ and $P \in S$ is such that $S \setminus P$ is not blocking, then there exists at least $2q + 1 - |S|$ lines through $P$ that pass through no other points of $S$. 
    \begin{proof}
        Suppose $m$ lines through $P$ pass through no other points of $S$. We may remove $P$ and add one point that is not $P$ on all but one of these lines to form a blocking set of $\AG(2, q)$. By Lemma \ref{lem:affine-blset}, this implies that $|S| + m - 2 \geq 2q - 1$ from which the result follows.
    \end{proof}
\end{lemma}

\begin{lemma}[Proposition 1.5 of \cite{except-one}]\label{lem:diamond-missing}
    Given $S \subseteq \Pp^2(\Fq)$ such that $|S| \leq 2q - 2$, if there exists no blocking set $S'$ of $\PG(2, q)$ such that $S' \supseteq S$ and $|S'| \leq |S| + 1$, then at least $2(2q - 1 - |S|)$ lines do not pass through any points of $S$.

    \begin{proof}
         Once again the goal is to construct a set containing $S$ that is incident to all lines except one. We use a greedy algorithm. Suppose initially that there are $m$ lines not passing through any point of $S$. At each step, if the lines not incident to the current set of points do not pass through a common point, add a point to the current set which is on a maximal number of those lines not already incident to the current set. Each step of this process strictly decreases the number of lines not incident to the constructed set, but can never decrease that number to $0$. Therefore, we eventually stop and obtain a situation where are lines not incident to the constructed set pass through a common point. 
         
         Say that up to this point we have added $k$ points to $S$, and there are still $x \geq 1$ lines remaining intersecting at a common point that are not incident. Then we achieve our desired set of points by adding for each of the $x$ lines but one, a point on the line that is not the common point to our constructed set. By our greediness, the fact that we never added the common point to our set implies that each of the first $k$ points added were on at least $x$ lines that were not already incident to our set. In fact, even if $x < 2$, each of them must have been on at least two, since the intersection point of any pair of lines always exists. We may therefore lower bound the number of lines not incident to $S$ originally as $m \geq \max(2, x)k + x$. By the hypothesis about the non-existence of a set $S'$, we have that $k \geq 1$ and therefore \[m - 2k - x \geq \max(0, x - 2)k \geq \max(0, x-2) \geq x - 2.\] Therefore, the number of points added in total to $S$, $k + x - 1$, is at most $\frac{m}{2}$. Hence, $|S| + \frac{m}{2} \geq |S| + k + x - 1 \geq 2q - 1$ via Lemma \ref{lem:affine-blset}, which establishes the result.
    \end{proof}
\end{lemma}

It turns out that the weaker bound obtained in Lemma \ref{lem:critical} is not strong enough on its own and it is therefore necessary to utilize Lemma \ref{lem:diamond-missing} for the majority of sets $S \subseteq \Pp^2(\Fq)$. As such, we define the following classes according to how close they are to being small blocking sets: \begin{align*}\spadesuit &\coloneqq \{S\subseteq \Pp^2(\Fq) \colon S \text{ blocking set such that } |S| < 2q\} \\ \clubsuit &\coloneqq \{S\subseteq \Pp^2(\Fq) \colon \exists\, T\in \spadesuit \text{ such that } |T \setminus S| = 1\} \\ \diamondsuit &\coloneqq \{S\subseteq \Pp^2(\Fq) \colon \nexists\, T\in \spadesuit \text{ such that } |T \setminus S| \leq 1\}.\end{align*} Then we have another summary result as a corollary:

\begin{corollary}\label{cor:M_S-bounds}
    \[|M_S| \geq \begin{cases} 2(2q - 1 - |S|) & \text{ if } S \in \diamondsuit \\ 2q - |S| & \text{ if } S \in \clubsuit \\ 0 & \text{ if } S \in \spadesuit \end{cases}.\]
\end{corollary}

We finally provide a proof of the new Theorem \ref{theorem:min-blset-bound}. It relies on an equivalent formulation of Lemma \ref{lem:affine-blset}.
\begin{corollary}\label{cor:reformat}
    Every blocking set, $S$, of size at most $2q$ contains a unique minimal blocking set.

    \begin{proof}
        If $S$ contains more than one, $S_1 \neq S_2$, then there is some line, $\ell$, that does not intersect $S_1 \cap S_2$ and therefore intersects both $S_1 \setminus S_2$ and $S_2 \setminus S_1$. But now $(S_1 \cup S_2) \setminus \ell$ is a set of less than $2q - 1$ points incident to all lines except $\ell$, contradicting Lemma \ref{lem:affine-blset}. We can moreover go in the reverse direction as well; given a set, $S$, of less than $2q - 1$ points incident to every line except $\ell$, choosing points $P \neq Q$ on $\ell$ makes $S \cup \{P, Q\}$ a set of size at most $2q$ that has minimal blocking sets containing $P$ but not $Q$ and vice-versa as subsets.
    \end{proof}
\end{corollary}
    \begin{proof}[Proof of Theorem \ref{theorem:min-blset-bound}]
        Our strategy is to exhibit two different ways of bounding the number of pairs $(B, S)$ in which $B$ is a minimal blocking set of size $k$ and $S \subseteq B$ is a subset of size $2k-2q$. Each minimal blocking set $B$ of size $k$ contains $\binom{k}{2k-2q}$ subsets of size $2k - 2q$, so there are exactly $B_k \binom{k}{2k-2q}$ such pairs. Now notice that every set of points of size $2k-2q$ must be included in at most one minimal blocking set of size $k$, as otherwise the union of two distinct such blocking sets would be a blocking set of size at most $2q$ containing both in contradiction to Corollary \ref{cor:reformat}. Therefore the number of such pairs is at most the number of sets of points of size $2k - 2q$, $\binom{q^2+q+1}{2k-2q}$. Combining these two approaches we see that $B_k \binom{k}{2k-2q} \leq \binom{q^2+q+1}{2k-2q}$ and obtain for $q \leq k \leq 2q$ the bound \[B_k \leq \frac{\binom{q^2 + q + 1}{2k - 2q}}{\binom{k}{2k-2q}}.\]
    \end{proof}

Notice that this bound, which we denote by $f(k)$, is actually sharp for $k = q + 1$. As $f$ is log-concave, making it unwieldy to bound directly later on, we will use the following unsharp version in the proof of the main theorem:

\begin{corollary}\label{cor:min-blset-bound}
    For $q+1 \leq k \leq 2q$, \[B_k < \left(\frac{e}{2}q\right)^{2k-2q}\]

    \begin{proof}
        We first check that $f(k)$ is log-concave in $k$ as already insinuated. Indeed, \[\frac{f(k+1)}{f(k)} = \frac{(q^2 + 3q + 1 - 2k)(q^2 + 3q - 2k)}{(2q-k)(k + 1)},\] is an increasing function for $q + 1 \leq k \leq 2q-1$. So on that same interval, any exponential function in $k$ is an upper bound over the interval if and only if it is an upper bound at the endpoints $\{q + 1, 2q\}$. Our claimed upper bound is also exponential, so we simply check the latter condition. The case of $k = q + 1$ can be verified easily using $q \geq 2$: \[f(q + 1) = q^2 + q + 1 \leq \frac{7}{4}q^2 \leq \frac{e^2}{4}q^2.\] For the other endpoint we use Stirling's Formula and the folklore fact that $(1 + x^{-1})^x < e$ for $x > 0$ to obtain that for $q \geq 3$, \begin{align*}f(2q) &= \binom{q^2 + q + 1}{2q} = \frac{1}{(2q)!}\prod_{k=0}^{q-1} (q^2 + q + 1 - k)(q^2 + q + 1 - (2q - 1 - k)) \\ &< \frac{[(q^2 + 1.5)^2]^q}{(2q)!} < \frac{1}{\sqrt{4\pi q}}\left(\frac{e(q^2 + 1.5)}{2q}\right)^{2q} < \frac{e^{\frac{3}{q}}}{\sqrt{4\pi q}}\left(\frac{e}{2}q\right)^{2q} < \frac{1}{\sqrt{q}}\left(\frac{e}{2}q\right)^{2q},\end{align*} which is more than sufficient. We check the $q = 2$ case separately, which reduces to the true statement $35 < e^4$.
    \end{proof}
\end{corollary}

\section{
Proof of Theorem \ref{theorem:TBS1}}\label{sec:main-proof}

We are now in a position to prove Theorem \ref{theorem:TBS1}. The main difficulty is in obtaining an upper bound on $\mu(\Omega)$. As suggested by the preliminary results in Section \ref{sec:Sam-Sung’s-Note}, the strategy is to divide and conquer based on the set of singularities of a curve $C_f$. To keep notation consistent, we will define the sets \begin{align*}\Sing_S &\coloneqq \{f \in R_\homog | \text{The set of rational singular points of $C_f$ is $S$}\} \\ \Sing'_S &\coloneqq \{f \in R_\homog | \text{$C_f$ is singular at every point in $S$}\} \\ \Sing'_{S, P} &\coloneqq \{f \in R_\homog | \text{$C_f$ is singular at every point in $S$, but not at $P$}\},\end{align*} the first and last of which were already dealt with in the previous section, and the second of which was considered informally in Section \ref{sec:TTFC}. Our bounding strategy is then encapsulated by the following covering: \[\Omega \subseteq \left(\bigcup_{S \in \diamondsuit} \Omega \cap \Sing_S \right) \cup \left(\bigcup_{S \in \spadesuit} \bigcup_{P \in S} \Omega \cap \Sing'_{S \setminus \{P\}, P} \right) \cup \left(\bigcup_{S \in \spadesuit} \Omega \cap \Sing'_S \right).\] Indeed, the set of rational singular points of a curve is either in $\diamondsuit$, or contains either all or all but one of a minimal blocking set of size at most $2q-1$.

\begin{theorem}\label{thm:upper-bound}
    Let $\mathbf{1}_\Z$ be the indicator function for the set of integers. Then, $\mu(\Omega) \leq (q^2+q+1)(q^{-3})^{q+1}\left(1 + e\left(\frac{e^2}{4q}\right)^{q^{\frac{2}{3}}}\right)+ \mathbf{1}_\Z(q^{\frac{1}{2}})(q^3 + q^{\frac{3}{2}})(q + 1)(q^{-3})^{q + q^{\frac{1}{2}} + 1} + e^{q+2}(q^{-1})^{4q-2}.$ 
    \begin{proof}
        The proof contains no new ideas; it is just a matter of simply unraveling the bounds already established and excessive bounding. Using the aforementioned covering we have \[\mu(\Omega) \leq \sum_{S \in \diamondsuit} \mu(\Omega \cap \Sing_S) + \sum_{S \in \spadesuit} \sum_{P \in S} \mu(\Omega \cap \Sing'_{S \setminus \{P\}, P}) + \sum_{S \in \spadesuit} \mu(\Omega \cap \Sing'_S).\] Now by applying Lemma \ref{lem:diamond-prob} and Lemma \ref{lem:club-prob} from Section \ref{sec:Sam-Sung’s-Note}, and using the bounds obtained from Lemma \ref{lem:tiny-B_k}, Corollary \ref{cor:M_S-bounds}, and Corollary \ref{cor:min-blset-bound}, we can bound each term separately. So we first have:
        \begin{align*}
            &\sum_{S \in \diamondsuit} \mu(\Omega \cap \Sing_S) \leq \sum_{S \in \diamondsuit} \left(\frac{1}{q}\right)^{3|S| + |M_S|} \leq \sum_{S \in \diamondsuit} \left(\frac{1}{q}\right)^{|S| + 4q - 2} \leq \left(\frac{1}{q}\right)^{4q - 2}\sum_{k=0}^{q^2+q+1}\sum_{\substack{S \in \diamondsuit \\ |S| = k}} \left(\frac{1}{q}\right)^k \\ &\leq \left(\frac{1}{q}\right)^{4q - 2}\sum_{k=0}^{q^2+q+1}\binom{q^2+q+1}{k} \left(\frac{1}{q}\right)^k \leq \left(\frac{1}{q}\right)^{4q - 2}\left(1 + \frac{1}{q}\right)^{q^2+q+1} \leq \left(\frac{1}{q}\right)^{4q - 2}e^{q+1}.
        \end{align*} For the next term we bound the final resulting geometric series by verifying by hand that it as most $e$ for $q \leq 5$, and that the result of the infinite sum, $(1 - \frac{4}{e^2}\cdot\frac{q+1}{q})^{-1}$ is at most $e$ for $q \geq 7$. Along the way we use that $q \geq 2$ and that $\frac{e}{2} < e^{\frac{1}{3}}$. In the last step we use $(e-1)e^{\frac{1}{3}q} \geq \left(\frac{3}{e}\right)e^{\frac{1}{3}q} = 3e^{\frac{1}{3}q - 1} \geq 3 \cdot \frac{q}{3} = q$.
        \begin{align*}&\quad\sum_{S \in \spadesuit} \sum_{P \in S} \mu(\Omega \cap \Sing'_{S \setminus \{P\}, P}) \leq \sum_{S \in \spadesuit} \sum_{P \in S} \left(\frac{1}{q}\right)^{3|S \setminus \{P\}|}\left(\frac{q+1}{q^2}\right)^{|M_{S \setminus \{P\}}|} \\&\leq \sum_{S \in \spadesuit} \sum_{P \in S} \left(\frac{1}{q}\right)^{3(|S|-1)}\left(\frac{q+1}{q^2}\right)^{2q + 1 - |S|} \\&\leq \left(\frac{q+1}{q^2}\right)^2\left(\frac{1}{q^3}\right)^{2q-2}\sum_{k=0}^{2q-1}\sum_{\substack{S \in \spadesuit \\ |S| = k}} k(q^2+q)^{2q-1-k}
        \end{align*}

        \begin{align*}
         &\leq (2q-1)\frac{(q+1)^2}{q^{6q-2}}\sum_{k=q+1}^{2q-1}\left(\frac{e}{2}q\right)^{2k-2q}(q^2+q)^{2q-1-k} \\ &\leq (2q-1)\frac{(q+1)^2}{q^{6q-2}}\left(\frac{e}{2}q\right)^{2q-2}\sum_{k=0}^{q-2}\left(\frac{4}{e^2} \cdot \frac{q+1}{q}\right)^k \leq (2q-1)\frac{(q+1)^2}{q^{6q-2}}\left(\frac{e}{2}q\right)^{2q-2}e \\ &\leq (4q+6)\left(\frac{1}{q}\right)^{4q-2}\left(\frac{e}{2}\right)^{2q-1} \leq 3\left(\frac{1}{q}\right)^{4q-3}\left(\frac{e}{2}\right)^{2q+2} \leq \left(\frac{1}{q}\right)^{4q-3}e^{\frac{2}{3}q+1} \\ &\leq (e-1)\left(\frac{1}{q}\right)^{4q-2}e^{q+1}.
        \end{align*} Similarly, to bound the final resulting geometric series in the following manipulation, note that it as at most $1$ for $q \leq 6$, and otherwise it is bounded by the infinite sum which is $(1 - \frac{e^2}{4q})^{-1}$ and at most $4e^{-1}$ for $q \geq 7$. 
        \begin{align*}
            &\sum_{S \in \spadesuit} \mu(\Omega \cap \Sing'_S) \leq \sum_{S \in \spadesuit} \left(\frac{1}{q}\right)^{3|S|} \leq \sum_{k=0}^{2q-1}\sum_{\substack{S \in \spadesuit \\ |S| = k}} \left(\frac{1}{q}\right)^{3k} \leq \left(\frac{1}{q}\right)^{3q}\sum_{k=0}^{2q-1} B_k\left(\frac{1}{q}\right)^{3k-3q} \\ &\leq B_{q+1}\left(\frac{1}{q^3}\right)^{q+1} + B_{q+q^{\frac{1}{2}}+1}\left(\frac{1}{q^3}\right)^{q+q^{\frac{1}{2}} + 1} + \sum_{k = \lceil q + q^{\frac{2}{3}} + 1 \rceil}^{2q-1}B_k\left(\frac{1}{q}\right)^{3k} \\ &\leq B_{q+1}\left(\frac{1}{q^3}\right)^{q+1}\left(1 + \frac{1}{q^2}\sum_{k = \lceil q + q^{\frac{2}{3}} + 1\rceil}^{q-2} \left(\frac{e}{2}q\right)^{2k-2q}\left(\frac{1}{q}\right)^{3(k-q-1)}\right) + B_{q+q^{\frac{1}{2}}+1}\left(\frac{1}{q^3}\right)^{q+q^{\frac{1}{2}}+1} \\ &\leq B_{q+1}\left(\frac{1}{q^3}\right)^{q+1}\left(1 + \left(\frac{e^2}{4q}\right)^{\lceil q^\frac{2}{3}\rceil}\cdot \frac{e^2}{4}\sum_{k=0}^{q - \lceil q^\frac{2}{3}\rceil - 2}\left(\frac{e^2}{4q}\right)^k\right) + B_{q+q^{\frac{1}{2}}+1}\left(\frac{1}{q^3}\right)^{q+q^{\frac{1}{2}}+1} \\ &\leq B_{q+1}\left(\frac{1}{q^3}\right)^{q+1}\left(1 + e\left(\frac{e^2}{4q}\right)^{\lceil q^{\frac{2}{3}} \rceil}\right) + B_{q+q^{\frac{1}{2}}+1}\left(\frac{1}{q^3}\right)^{q+q^{\frac{1}{2}}+1} \\ &\leq (q^2+q+1)\left(\frac{1}{q^3}\right)^{q+1}\left(1 + e\left(\frac{e^2}{4q}\right)^{q^{\frac{2}{3}}}\right) + \mathbf{1}_\Z(q^{\frac{1}{2}})(q^3 + q^{\frac{3}{2}})(q + 1)\left(\frac{1}{q^3}\right)^{q+q^{\frac{1}{2}}+1}.
        \end{align*}
         Now summing gives the bound we claimed.
    \end{proof}
\end{theorem}

On the other hand, the prominent first two terms of this upper bound correspond to readily exhibited classes of transverse-free curves. Applying the exact same strategy as in Theorem \ref{thm:lubber}, we obtain the following lower bound:

\begin{theorem}\label{thm:lower_bound}
    $\mu(\Omega) \geq (q^2+q+1)q^{-3(q+1)}(1-q^{-3q+2}) + \mathbf{1}_\Z(q^{\frac{1}{2}})(q^3 + q^{\frac{3}{2}})(q + 1)q^{-3(q+q^{\frac{1}{2}}+1)}(1 - q^{-3q+3q^{\frac{1}{2}}+4}).$

    \begin{proof}
        We will use the inclusion-exclusion principle to lower bound the density of the set of curves for which the rational points of some line or Baer subplane are singularities. As a consequence of Corollary \ref{cor:reformat}, the union of any two distinct blocking sets must contain at least $2q + 1$ points. Now by using this fact and truncating after the second set of terms, our bound is \begin{align*}& B_{q+1}q^{-3(q+1)} + \mathbf{1}_\Z(q^{\frac{1}{2}})B_{q + q^\frac{1}{2} + 1}q^{-3(q+q^{\frac{1}{2}}+1)} - \binom{B_{q+1}+ \mathbf{1}_\Z(q^{\frac{1}{2}})B_{q + q^\frac{1}{2} + 1}}{2}q^{-3(2q+1)} \\ = \,\,& (q^2+q+1)q^{-3(q+1)} - \binom{q^2+q+1}{2}q^{-3(2q+1)} + \mathbf{1}_\Z(q^{\frac{1}{2}})(q^3 + q^{\frac{3}{2}})(q + 1)q^{-3(q+q^{\frac{1}{2}}+1)} \\ &\qquad- \mathbf{1}_\Z(q^{\frac{1}{2}}) [(q^3 + q^{\frac{3}{2}})(q + 1)(q^2+q+\tfrac{1}{2}+ \tfrac{1}{2}(q^3 + q^{\frac{3}{2}})(q + 1))q^{-3(2q+1)}].\end{align*} The claimed bound may then be obtained by using $\frac{1}{2}(q^2 + q) \leq q^2$ for $q \geq 2$ and $q^2+q+\frac{1}{2}+ \frac{1}{2}(q^3 + q^{\frac{3}{2}})(q + 1) \leq q^4$ for $q \geq 2^2$.
    \end{proof}
\end{theorem}

Combining the upper and lower bounds on $\mu(\Omega)$ just obtained, we see that:

\begin{corollary}\label{cor:gen-asympt-bound}
    \[\mu(\Omega) = ((q^2+q+1)(q^{-3})^{q+1} + \mathbf{1}_{\Z}(q^{\frac{1}{2}})(q^3 + q^{\frac{3}{2}})(q + 1)(q^{-3})^{q + q^{\frac{1}{2}} + 1})\left(1+O\left(\left(4e^{-2}q\right)^{-q^{\frac{2}{3}}}\right)\right).\]
\end{corollary}

\begin{proof}[Proof of Theorem \ref{theorem:TBS1}]
    Since $4e^{-2} > \frac{1}{2}$, Corollary \ref{cor:gen-asympt-bound} gives us the claimed asymptotic formula for general $q$. For the particular cases of $q = p$ and $q = p^2$, by applying the result from \cite{best-bounds}, we have that the exponent $e$ is precisely either $k$ or $k/2$ and so all blocking sets of size less than $\frac{3}{2}(q+1)$ must either be the $\Fq$ points of a line or a Baer subplane. Hence $B_k = 0$ for all $k < \frac{3}{2}(q+1)$ that are not $q+1$ or $q+q^{\frac{1}{2}}+1$ instead of the general bound of $k < q + q^{\frac{2}{3}} + 1$, and re-running the proof of Theorem \ref{thm:upper-bound} obtains the same result but with $q^{\frac{2}{3}}$ replaced by $\frac{1}{2}(q+1)$. Similarly, Corollary \ref{cor:gen-asympt-bound} is modified by with the same replacement of $q^{\frac{2}{3}}$ with $\frac{1}{2}(q+1)$. Then using $4e^{-2} > \frac{1}{2}$ again and simplifying $\mathbf{1}_{\Z}(q^{\frac{1}{2}})$ gives the claimed asymptotic formula in these special cases.
\end{proof}

\begin{remark}
    From the efficiency heuristic in Section \ref{sec:TTFC}, we strongly believe that the claim expressed in the title of this paper generalizes to higher dimensions. More precisely, fix positive integers $k \leq n$. Then for sufficiently large $q$, almost every hypersurface in $\Pp^n$ that is not transverse to all embedded $\Pp^k$ should be singular at every rational point of an embedded $\Pp^{n-k}$ for $1 \leq k \leq n$. However we were not able to provide a proof of this claim since several of the key ingredients in our proofs seem specific to the two-dimensional case. As already mentioned, it does not appear that Theorem \ref{thm:inQ} of Subsection \ref{subsec: rationality} generalizes, so we may be required to work with the densities in an approximate matter. In addition, we do not know of any helpful generalization of Theorem \ref{thm:Sam-Sung’s-Note} to greatly reduce casework. Crucially however, there does not exist sufficient results in the literature on the number of ways of covering projective spaces with linear subspaces that would generalize Section \ref{sec:block-time}. 
\end{remark}

\bibliographystyle{alpha}
\footnotesize{\bibliography{refs}}

\end{document}